\documentclass[3p,10pt]{elsarticle}

\usepackage{graphicx}
\usepackage{color}

\usepackage{epsfig,amssymb}
\usepackage{amsmath}
\usepackage{capt-of}

\usepackage{multicol}
\usepackage{multirow}

\def\doubleone{ 1 \kern-.224em \hbox{\rm l}}
\def\iu{{\rm i}}
\def\e{{\rm e}}
\def\d{{\rm d}}

\def\doubleone{ 1 \kern-.224em \hbox{\rm l}}
\def\iu{{\rm i}}
\def\e{{\rm e}}
\def\d{{\rm d}}

\newtheorem{theorem}{Theorem}
\newtheorem{lemma}{Lemma}

\newtheorem{remark}{Remark}

\newcommand{\norm}[1]{\left\Vert#1\right\Vert}

\begin{document}

\renewcommand\arraystretch{1.0}

\begin{frontmatter}

\title
{Convergence in the maximum norm of ADI-type methods for parabolic problems}

\author[Lag]{S. Gonz\'{a}lez-Pinto and  D.
Hern\'{a}ndez-Abreu \fnref{aut1} }

\fntext[aut1]{This work has been partially supported by the Spanish Project
MTM2016-77735-C3-3-P of  Ministerio de Econom\'{i}a,
Industria y Competitividad.}

%\cortext[cor1]{Corresponding author: D. Hern\'{a}ndez
%Abreu}

\baselineskip=0.9\normalbaselineskip
\vspace{-3pt}

%\maketitle

\address[Lag]{ { \footnotesize Departamento de An\'{a}lisis Matem\'{a}tico. Universidad de La
Laguna. 38071. La Laguna, Spain. \\ email:
spinto\symbol{'100}ull.edu.es, dhabreu\symbol{'100}ull.edu.es}}

%\maketitle

\begin{abstract}
Results on unconditional convergence in the Maximum norm for ADI-type methods, such as the Douglas method, applied to the time integration of semilinear parabolic problems are quite difficult to get, mainly when the number of space dimensions $m$ is greater than two. Such a result is obtained here under quite general conditions on the PDE problem in case that  time-independent Dirichlet boundary conditions are imposed. To get these bounds, a theorem that guarantees, in some sense,  power-boundeness of the stability function  independently of both the space and time resolutions is proved. 
\end{abstract}


\end{frontmatter}

\section{Introduction}
\noindent The present article considers the numerical solution of ODE systems
\begin{equation}\label{ODEm}
\dot U = D\, U + g(t) , \quad U(0)=U_0, \quad  t\in [0,t^*],\qquad D:=\sum_{j=1}^m D_j, \quad g(t):=\sum_{j=1}^m g_j(t),\end{equation}
stemming from the spatial discretization by using finite differences (or finite volumes) of semilinear parabolic PDEs with constant diffusion coefficients $\beta_j>0$  and an initial condition and  Dirichlet Boundary Conditions (BCs) of the form
\begin{equation}\label{lindiff}\begin{array}{c}
\partial_t u(t,\vec{x}) = \displaystyle \sum_{j=1}^m \beta_j\,\partial_{x_j x_j} u(t,\vec{x}) +   c(t,\vec{x}), 
\quad t\in [0,t^*] ,\quad\vec{x}=(x_1,\ldots,x_m)^\top\in \mathcal{I}:=(0,1)^m,\\ 
u(0, \vec{x})=u_0(\vec{x}), \; \vec{x}\in \partial \mathcal{I}; \qquad u(t,\vec{x})=h(t),\; t \in (0,t^*],\; \vec{x}\in \partial \mathcal{I}. \end{array} \end{equation}
Here, $c(t,\vec{x})$ is a source term and we assume that its discretization is entirely included in $g_1(t)$, so that $g_j(t), \:j=2,\ldots,m$, consist only of contributions from the boundary conditions in the $j$-direction. In particular, we shall be concerned with time-independent Dirichlet boundary conditions, in which case the vectors $g_2(t),\ldots , g_m(t)$ are constant.

\smallskip

\noindent To prove convergence  in the maximum norm for many numerical methods of splitting type applied to (\ref{ODEm}), it is customary to get uniform bounds  for $\|R(\tau D_1,\ldots,\tau D_m)^n\|_\infty, \; n=1,2,\ldots,t^*/\tau$, where $\tau>0$ is the time stepsize and $R(\tau D_1, \ldots,\tau D_m)$ is a rational mapping acting on the matrices $\tau D_j$, $1\leq j\leq m$. Typically,  we have $D_j= \beta_j\, (I_{N_{m}}\otimes\ldots\otimes L_j\otimes\ldots \otimes I_{N_{1}})$, with $ L_j=\mbox{\rm tridiag}(1,-2,1)/\Delta  x_j^2$ when second order central differences are considered in the spatial discretization of (\ref{lindiff}).  Here, we denote the spacing $\Delta x_j = 1/(N_{j}+1)$, where  $N_j$ is the number of equidistant grid-points on the $j$-direction, and $\otimes$ stands for the Kronecker product of matrices.  It should be observed that the matrices $D_j$ pairwise commute.  Such methods of splitting type when applied to (\ref{ODEm}) typically produce a recursion for the global errors $E_n:=U_n-U(t_n)$, $n\geq 0$, of the form $E_{n+1}=RE_n+S_n$, where $U_n$ stands for the numerical solution at $t_n=n\tau $, $S_n$ denotes the local error and $R$ is the  stability matrix associated to the numerical integrator. For ADI-type integrators the stability matrix depends on $D_1,\ldots, D_m$  (see, e.g., \cite[Sec. II.2.3]{Hun-Ver03}). A relevant example is
\begin{equation}\label{stab-matrix-m}
R(\tau D_1, \ldots,\tau D_m ) = I + \Pi (\theta )^{-1} \tau D , \quad D=\sum_{j=1}^m D_j,
\end{equation}
where $ \Pi (\theta ) = (I - \theta \tau D_1)\cdots (I-\theta \tau D_m )$,  which has the associated 
stability function of $m$ complex variables
\begin{equation}\label{stab-functm}
R(z_1, z_2,\ldots,z_m)  = 1  + \frac{z}{Q}, \quad z:=\sum_{j=1}^m z_j, \quad Q:=\prod_{j=1}^m (1-\theta z_j).
\end{equation}  For the choice
$\theta = 1/2$, this is the stability matrix of  the Peaceman-Rachford method (when $m=2$), also   the one of the  Douglas scheme (\cite{Douglas62}, \cite{Hun-98a},~\cite[p. 373] {Hun-Ver03})
and  the  one of the one-stage AMF-W-method
\cite{GHH-sinum20}. Furthermore, the stability matrix of the so-called
Hundsdorfer--Verwer scheme \cite[Section~IV.5.2]{Hun-Ver03},
which is a $2$-stage W-method
of order $2$ in general, and of order $3$ for $\theta = (3+\sqrt 3)/6$,
is given by
%\begin{equation}\label{hd-m}
$$R(\tau D_1, \ldots,\tau D_m ) = I + 2 \Pi (\theta )^{-1} \tau D - \Pi (\theta )^{-2} \tau D 
+ \frac 12 \bigl( \Pi (\theta )^{-1} \tau D \bigr)^2 .$$
%\end{equation}
In this case the stability function is given by  
%\begin{equation}\label{RHV}
$R(z_1, z_2,\ldots,z_m)  = 1  + \frac{2z}{Q} + \frac{z^2-2z}{2Q^2}$
%\end{equation}
(with $z$ and $Q$ defined in (\ref{stab-functm})). The power boundedness in the maximum norm of some $R$-stability functions was already considered in \cite{GHH-BIT21}. However,  for $m\ge 3$ the power bound  there obtained is not uniform, since it allows a logarithmic growth in terms of $n$ or $N$ \cite{GHH-BIT21}, i.e.,
$$ \|R(\tau D_1, \ldots,\tau D_m)^n\|_\infty \le C_m \min\{(\ln(n+1))^m,\:(\ln(N+1))^m\}, \; \tau>0,\; n\tau \le t^*, \; \;N=\max\{N_j,\:j=1,\ldots,m\}.$$ 
With this power bound,  convergence results in the maximum norm of size  $\mathcal{O}(\tau^p|\log(\tau)|^m)$, when the local errors are of size $\mathcal{O}(\tau^{p+1})$ can be obtained. However, with power bounds of the stability matrix as the one in (\ref{4-eq1}), it can be shown    convergence  of size  $\mathcal{O}(\tau^p)$ in case of  time independent BCs in (\ref{lindiff}). In Section \ref{sect2}, we prove a result related to the power boundedness for rational functions. This result is applied in Section \ref{sect3} to show unconditional convergence in the maximum norm for some ADI-methods. In Section 4,  numerical experiments are included to illustrate the orders of convergence regarding the PDE solution for some relevant ADI-type methods.

\section{Bounds in the maximum norm for rational functions}\label{sect2}
\noindent We look for bounds in the maximum norm of the form
\begin{equation}\label{4-eq1}\begin{array}{c}
\|D^{-1} R(\tau D_1, \tau D_2, \ldots, \tau D_m)^n\|_\infty \le C < \infty, \quad t_n=n \tau, \;n=0,1,2,\ldots, \;\tau>0, \\[0.5pc]   D=\displaystyle\sum_{j=1}^m D_j, \quad  D_j=\beta_j(I_m\otimes \cdots \otimes I_{j+1} \otimes L_j\otimes  I_{j-1} \otimes \cdots \otimes I_1), \\[0.5pc] L_j= \frac{1}{\Delta x_j^2} \hbox{\rm tridiag}(1,\:-2,\:1)_{N_j}, \quad I_j \;\hbox{\rm is the identity matrix of dimension } N_j,\quad  \Delta x_j=\frac{1}{N_j+1},
\end{array}
\end{equation}
where $\displaystyle R(z_1,\ldots,z_m)$ is a rational function (or a mapping when acting on the matrices $D_j$) of $m$ complex variables that is $A_m(\alpha)$-stable, i.e.
\begin{equation}\label{4-eq2}\begin{array}{c}
|R(z_1,\ldots,z_m)| \le 1, \quad   \forall z_j \in \mathcal{W}(\alpha):=\{z\in \mathbb{C}:\: |\arg (-z)| \le \alpha\}\cup\{0\},\quad \hbox{\rm for some } 0<\alpha \le \pi/4.
% \\ Q(z_1,\ldots,z_m)=0 \Longrightarrow z_j \in K,   \;\forall j=1,\ldots,m.
\end{array}
\end{equation} 
Of course, if (\ref{4-eq2}) holds true for some $\bar{\alpha}\in (0,\pi/2]$ then it also holds for $\alpha=\min\{\bar{\alpha},\:\pi/4\}$.  
\begin{theorem}\label{4-eq3}
If $R(z_1,z_2,\ldots,z_m)$ is a rational function  that satisfies (\ref{4-eq2}), then there exists a constant $K$ only depending on   $\alpha$ and $m$ such that  (\ref{4-eq1}) holds for $C=K/\hat{\beta}, \:\hat{\beta}=:\min_{1\le j \le m} \beta_j>0$. 
\end{theorem}  
The proof of this theorem is given below and makes use of the following two lemmas. 
\begin{lemma}\label{4-eq4}
For any matrix $L_j$ (\ref{4-eq1}) it holds that 
\begin{equation}\label{4-eq5}\| L_j^{-1}\|_\infty \le \frac{1}{8}, \end{equation} 
\begin{equation}\label{4-eq5a}\| (zI_j- \tau L_j)^{-1}\|_\infty \le \frac{\sec (\theta/2) }{|z|}, \quad \forall \tau>0 \;\hbox{\rm and } z=|z| e^{\iu \theta}, \;|\theta| < \pi, \end{equation}   
\begin{equation}\label{4-eq6} \| (zI_j - \tau L_j)^{-1}\|_\infty \le \frac{1}{8\tau - |z|}, \quad \forall |z| <8\tau,
\end{equation} 
and 
\begin{equation}\label{4-eq6a} \| (zI_j - \tau L_j)^{-1}\|_\infty \le \frac{1}{ |z|- 4 \tau/\Delta x_j^2}, \quad \forall |z| > 4\tau/\Delta x_j^2>0.
\end{equation} 

\end{lemma}
{\bf Proof.} The formula in (\ref{4-eq5}) is well known in the literature (see, e.g., \cite[formula (4.10)]{Mattheij} or \cite[p. 43-45]{lar-thom}). The formula (\ref{4-eq5a}) is  an immediate consequence of Lemma 4.1 in  \cite{GHH-BIT21} (see also \cite[formula (5)]{farago02ste}), with $\mu=\tau/(\Delta x_j)^2$ and $T_j=\mbox{\rm tridiag}(1,-2,1)_{N_j}$, since 
$$  \| (zI_j- \tau L_j)^{-1}\|_\infty  = \mu^{-1}\| (\mu^{-1}z I_j-  T_j)^{-1}\|_\infty \le  \mu^{-1} \frac{\sec (\arg (\mu^{-1} z) /2) }{\mu^{-1}|z|} = \frac{\sec (\theta/2) }{|z|}.$$ 
To show (\ref{4-eq6}), by considering $ |z| <8\tau$, we have that 
$$\displaystyle  \| (zI_j - \tau L_j)^{-1}\|_\infty = \| (\tau L_j)^{-1} (I_j - z\tau^{-1} L_j^{-1})^{-1}\|_\infty \le \frac{\tau^{-1}\|L_j^{-1}\|_\infty }{1 - |z| \tau^{-1} \|L_j^{-1}\|_\infty}\le \frac{\tau^{-1}8^{-1}}{1 - |z| \tau^{-1} 8^{-1}} = \frac{1}{8\tau - |z|}.$$
To show (\ref{4-eq6a}), for $|z| > 4\tau/\Delta x_j^2$ it holds
$$\displaystyle  \| (zI_j - \tau L_j)^{-1}\|_\infty = \| z^{-1} (I_j - z^{-1}\tau L_j)^{-1}\|_\infty \le \frac{|z|^{-1} }{1 - |z|^{-1} \tau \|L_j\|_\infty}\le \frac{1}{  |z|- 4\tau /\Delta x_j^2}.$$ 

\hfill $\Box$ 

\begin{lemma}\label{4-eq15}
Assume that for any positive integer $m$ we have that  
\begin{equation}\label{complexnumbers}
z_j=-r_je^{\iu\theta_j  }, \quad r_j>0,\;-\alpha \le \theta_j \le \alpha,\;(j=1,2,\ldots,m), \quad 0\le \alpha \le \pi/4.
\end{equation}
Then 
\begin{equation}\label{inequalities}
\displaystyle |\sum_{j=1}^m z_j  | \ge \sqrt{\sum_{j=1}^m |z_j|^2} \ge \frac{1}{\sqrt{m}}\sum_{j=1}^m |z_j| \ge \sqrt{m} \prod_{j=1}^m |z_j|^{1/m}.
\end{equation}
\end{lemma}
{\bf Proof.} The last two inequalities in (\ref{inequalities}) follow from the fact that for positive numbers the Quadratic Mean is greater or equal  than the Arithmetic Mean and this is greater or equal than the Geometric Mean. To show the first inequality, we observe that for  complex numbers $z_1,z_2$ satisfying (\ref{complexnumbers}) it holds that
$|z_1+z_2|^2=|z_1|^2+|z_2|^2+2{\rm Re}(z_1\overline{z}_2)\geq |z_1|^2+|z_2|^2.$
Hence, $s_2:=z_1+z_2$ fulfils $|s_2|\geq \sqrt{|z_1|^2+|z_2|^2}$ and it has an angle $|\hat{\theta}_2|\leq \alpha$ $(\leq \pi/4)$ with the negative $x-$axis. In particular, $s_2$ takes the form (\ref{complexnumbers}). Then, adding a new complex number $z_3$ (\ref{complexnumbers}) and using the same argument we deduce that $s_3:=s_2+z_3$ fulfils $|s_3|\geq \sqrt{|s_2|^2+|z_3|^2}\geq \sqrt{|z_1|^2+|z_2|^2+|s_3|^2}$ and it has an angle $|\hat{\theta}_3|\leq \alpha$ $(\leq \pi/4)$ with the negative $x-$axis. The application of the induction principle concludes the proof. \hfill $\Box$

\medskip

\noindent {\bf Proof of Theorem \ref{4-eq3}.} We define $\beta:=\max\{\beta_1,\ldots,\beta_m\}$ and  use below the following notation for the Kronecker product of matrices 
$ \displaystyle{[\otimes A_l]_{l=1}^{k}} := A_{k}\otimes A_{k-1}\otimes \ldots \otimes A_1.$   Consider the positively oriented boundary   $\Gamma = \gamma_1 \cup \gamma_2 \cup \gamma_3$ of the open domain $\Omega\subset \{z:\;\hbox{Re } z \le 0\}$, which is  symmetric with respect to the negative real axis in the complex plane,
\begin{equation}\label{4-eq7}
\begin{array}{c}
\gamma_1 := \{ -r \e^{ \iu\theta }\, ; \, -\alpha\le \theta \le \alpha\},\quad \gamma_3:= \{- r^* \e^{\iu \theta }\, ; \, -\alpha\le \theta\le \alpha \},  \\[2mm]
\gamma_2 := \{ -\rho \e^{  - \iu \alpha  },\: r \le \rho \le r^* \} \cup  \{ -(r+r^*-\rho) \e^{   \iu \alpha },\: r \le \rho \le r^* \},  
\\[2mm]
 r:=4\hat{\beta} \tau ,\qquad   r^*:=8\beta\tau /(\Delta x)^2, \quad \Delta x :=\min\{\Delta x_1,\ldots,\Delta x_m\}<1, \;\hbox{\rm since each } \Delta x_j=1/(N_j+1)< 1.
\end{array}
\end{equation} 
Observe that $0<r< \frac{ r^*}{2}.$ 
Let us define the rational function (and the associated mapping when acting on matrices)
\begin{equation}\label{4-eq8}
\phi(z_1,z_2,\ldots,z_m) = (z_1+z_2+ \ldots +z_m)^{-1}R(z_1,z_2,\ldots,z_m)^n.
\end{equation}
Taking into account that $\phi(z_1,z_2,\ldots,z_m)$ is analytic if $ (z_1,z_2,\ldots,z_m)\in \bar \Omega^m$, from the Cauchy's integral formula applied on each variable we get the following formula by using iterated  integrals 
\begin{equation}\label{4-eq9}\phi(z_1^*,z_2^*,\ldots,z_m^*)=\frac{1}{(2\pi \iu)^m} \oint_\Gamma\cdots \oint_\Gamma  \phi(z_1,z_2,\ldots,z_m) \prod_{j=1}^m\frac{dz_j}{z_j-z_j^*}, \quad \forall (z_1^*,z_2^*,\ldots,z_m^*) \in \Omega^m.
\end{equation}
By considering the mapping acting on the matrices $D_j$ we deduce that  
\begin{equation}\label{4-eq10}\phi(\tau D_1,\ldots,\tau D_m)=\frac{1}{(2\pi \iu)^m} \oint_\Gamma\cdots \oint_\Gamma  \phi(z_1,z_2,\ldots,z_m) \displaystyle{[\otimes (z_j I_j - \tau \beta_jL_j)^{-1} dz_j]_{j=1}^{m}}.
\end{equation} 
Observe that the eigenvalues of each matrix $\tau \beta_j L_j$ are $\lambda_i^{(j)}=-4\beta_j \frac{\tau}{(\Delta x_j)^2}\sin\left(\frac{\pi}{2}i\Delta x_j\right)^2$, $i=1,\ldots,N_j$, where  
$$ -r^*=-\frac{8\beta \tau}{(\Delta x)^2} < \lambda_{N_j}^{(j)} < \lambda_{1}^{(j)}=-4\beta_j \tau \left(\frac{\sin (\pi/2\cdot\Delta x_j)}{\Delta x_j}\right)^2<-4\beta_j \tau \le -r.  $$ 
Hence, the spectrum of $\tau \beta_j L_j$ falls  in  $(-r^*,-r)\subset\Omega$.

\smallskip

\noindent At this point we should notice  the identity  $D^{-1} R(\tau D_1, \tau D_2, \ldots, \tau D_m)^n= \tau \phi(\tau D_1,\ldots,\tau D_m).$ From here, taking the maximum norm and using that $\norm{M_1\otimes M_2}_\infty =\norm{M_1}_\infty\norm{M_2}_\infty$ for two matrices $M_1$ and $M_2$, we get that 
\begin{equation}\label{4-eq10a}\begin{array}{l}
\|D^{-1} R(\tau D_1, \tau D_2, \ldots, \tau D_m)^n\|_\infty \le  \displaystyle {\bf A}:=
 \frac{\tau}{(2\pi)^m} \oint_\Gamma\cdots \oint_\Gamma  |\phi(z_1,z_2,\ldots,z_m)| \prod_{j=1}^{m} \|(z_j I_j - \tau \beta_j L_j)^{-1}\|_\infty |dz_j|. \end{array}
\end{equation} 

\noindent Next we bound  $\|(z_j I_j - \tau \beta_j L_j)^{-1} dz_j\|_\infty$ when $z_j \in \Gamma$. We distinguish three cases. 
\begin{enumerate}
\item $z_j\in \gamma_1$, then $z_j=-r e^{\iu \theta_j}, \; -\alpha \le \theta_j \le \alpha$, and $|dz_j|=r d\theta_j$. From (\ref{4-eq6}) it follows that
$$\displaystyle  \|(z_j I_j - \tau\beta_j  L_j)^{-1}\|_\infty \le \frac{1}{8\tau\beta_j - r}\le    \frac{1}{8\tau\hat{\beta} - r} =\frac{1}{r}$$ and we deduce that
\begin{equation}\label{4-eq12}
z_j\in \gamma_1 \Longrightarrow  \|(z_j I_j - \tau \beta_j L_j)^{-1}\|_\infty |dz_j| \le \frac{r}{r}d\theta_j=d \theta_j.
\end{equation}  
\item $z_j\in \gamma_2$, then $z_j=-\rho_j e^{ \pm \iu\alpha_j}, \; r\le \rho_j \le r^*$, and $|dz_j|=d\rho_j$. From (\ref{4-eq5a}) it follows that 
$$\displaystyle  \|(z_j I_j - \tau\beta_j  L_j)^{-1}\|_\infty \le \frac{\sec{((\pi-\alpha)/2)}}{|z_j|}=   \frac{1}{\rho_j \sin (\alpha/2)}$$ and then
\begin{equation}\label{4-eq13}
z_j\in \gamma_2 \Longrightarrow  \|(z_j I_j - \tau \beta_j L_j)^{-1}\|_\infty |dz_j| \le \frac{1}{\rho_j \sin (\alpha/2)}d\rho_j.
\end{equation}    
\item $z_j\in \gamma_3$, then $z_j=-r^* e^{ \iu\theta_j}, \; -\alpha \le \theta_j \le \alpha$, and $|dz_j|=r^* d\theta_j$. From (\ref{4-eq6a}) it follows that
$$\displaystyle  \|(z_j I_j - \tau\beta_j L_j)^{-1}\|_\infty \le \frac{1}{|z_j| - 4\tau\beta_j /\Delta x_j^2} \le \frac{1}{r^* - r^*/2}=\frac{2}{r^*}$$ and we get that
\begin{equation}\label{4-eq14}
z_j\in \gamma_3 \Longrightarrow  \|(z_j I_j - \tau\beta_j L_j)^{-1}\|_\infty |dz_j| \le 2 d\theta_j.
\end{equation}  
\end{enumerate}
From the A($\alpha$)-stability of $R(z_1,\ldots,z_m)$  we  deduce that 
\begin{equation}\label{4-eq11}
|\phi(z_1,z_2,\ldots,z_m)| \le  \frac{1}{|z_1+ \ldots+ z_m|}, \quad \forall (z_1,\ldots,z_m) \in (\bar{\Omega})^m.  
\end{equation}  
From (\ref{inequalities}) in Lemma 2, we have that 
$ |z_1+ \ldots+ z_m|\ge \sqrt{m} \prod_{j=1}^m |z_j|^{1/m},   \quad \forall (z_1,\ldots,z_m) \in (\bar{\Omega})^m.$ Consequently, this together with  (\ref{4-eq11}) yields
\begin{equation}\label{4-eq16}
|\phi(z_1,z_2,\ldots,z_m)| \le m^{-1/2} \prod_{j=1}^m |z_j|^{-1/m}, \quad \forall (z_1,\ldots,z_m) \in (\bar{\Omega})^m.  
\end{equation}  
Now, by considering (\ref{4-eq10a})  and  (\ref{4-eq16}) we deduce that 
\begin{equation}\label{4-eq17} \begin{array}{c}\displaystyle \hbox{\bf A}\le\hbox{\bf B}:=  \frac{\tau }{\sqrt{m}(2\pi)^m}\sum_{1\le i_1,i_2,\ldots,i_m \le 3}\int_{\gamma_{i_1}} \cdots \int_{\gamma_{i_m}}   \prod_{j=1}^{m} \left( \|(z_j I_j - \tau \beta_j L_j)^{-1}\|_\infty |z_j|^{-1/m} |dz_j|\right) .
\end{array} 
\end{equation}  
Taking account that all these iterated integrals can be transformed into products of integrals in one variable, we get 
\begin{equation}\label{4-eq18} \begin{array}{l}\displaystyle \hbox{\bf B} \le \frac{\tau }{\sqrt{m}(2\pi)^m} \sum_{m_1 + m_2 +m_3=m\atop m_j\ge 0}  \frac{m!}{m_1! m_2!  m_3!} (A_1)^{m_1}(A_2)^{m_2}(A_3)^{m_3}, \quad \hbox{\rm where, using (\ref{4-eq12}), (\ref{4-eq13}) and (\ref{4-eq14}) respectively,}\\[1pc]
\displaystyle A_1 :=  \int_{-\alpha}^ \alpha r^{-1/m} d\theta =2\alpha r^{-1/m},
\\[1pc]
\displaystyle A_2:= \int_r^{r^*} \frac{\rho^{-1-1/m}}{\sin (\alpha/2)} d\rho = \frac{m}{\sin (\alpha/2)}\left(r^{-1/m} -  (r^*)^{-1/m}\right) < \frac{m }{\sin (\alpha/2)}r^{-1/m},\quad {\rm and}
 \\[1pc] \displaystyle A_3 := \int_{-\alpha}^\alpha (r^*)^{-1/m} 2 d\theta   =4\alpha (r^*)^{-1/m}. \end{array} 
\end{equation} 
Then, we have for $m_1+m_2+m_3=m$,
$$\displaystyle \tau  (A_1)^{m_1}(A_2)^{m_2}(A_3)^{m_3} \le  \tau\:C_{m_1,m_2,m_3}\: r^{-(m_1+m_2)/m}(r^*)^{-m_3/m},$$ with $$  \displaystyle  C_{m_1,m_2,m_3}:=\left(2 \alpha \right)^{m_1} \left(\frac{m}{\sin(\alpha/2)}\right)^{m_2} \left(4\alpha\right)^{m_3},
$$
and
$$
 r^{-(m_1+m_2)/m}(r^*)^{-m_3/m}=\displaystyle{ r^{-1} \left(\frac{r}{r^*}\right)^{m_3/m}}\le \displaystyle{ \frac{1}{4\hat{\beta}\tau }\left(\frac{1}{2}\right)^{m_3/m}}.$$ Hence, each term $\tau  (A_1)^{m_1}(A_2)^{m_2}(A_3)^{m_3}$ is bounded since
$$  \tau  (A_1)^{m_1}(A_2)^{m_2}(A_3)^{m_3} \le  \hat \beta^{-1} \frac{C_{m_1,m_2,m_3}}{4 \cdot 2^{m_3/m}}.$$
This concludes the proof. \hfill $\Box$

\section{Convergence in the uniform norm of some ADI-type methods}\label{sect3}

\noindent The first goal of this section is to show unconditional convergence of order two in the maximum norm for semilinear parabolic problems with constant diffusion coefficients (and a time dependent source term) and time-independent Dirichlet boundary conditions (\ref{ODEm})-(\ref{lindiff}), when the one-step AMF-W-method (henceforth denoted as {\bf AMF-W1}) in \cite{GHH-sinum20,gonzalez18amf} is considered  with the parameter choice $\theta=1/2$
\begin{equation}\label{AMF-W1-m}
\begin{aligned}
%\begin{array}{rcl}
K_1^{(0)} & = \displaystyle \tau D\,  U_n 
+ \tau\, g(t_n),   \\[1.mm]
(I-\theta \tau D_j ) K_1^{(j)} & = K_1^{(j-1)} + \theta \tau^2 \dot g_j (t_n),\quad j=1,\ldots,m,\\[1.mm]
\displaystyle U_{n+1} &  = U_n + K_1^{(m)},
%\end{array}
\end{aligned}
\end{equation}
where $\dot v(t)$ stands for the derivative of a function $v(t)$ regarding $t$. The following discussion can be applied in similar terms to the Douglas method \cite[p. 373] {Hun-Ver03}.
We use the same notations as in \cite{GHH-sinum20}. 
The global error  at the time step $t_n=n\tau $ is denoted as in \cite[formula (2.3)]{GHH-sinum20} by $ E_n=U_n-U(t_n),$
where $U_n$ is the solution of the numerical method and $U(t)=u(t,\vec{x}_G)$ is at the same time the exact solution of the  (\ref{ODEm}) and the exact solution of the PDE on the set of discrete points $\vec{x}_G$ of the spatial mesh-grid $G$. Observe that we will not consider in our analysis  the truncation errors introduced in the spatial discretization of the PDE, since when using central differences we get a stable space discretization and the truncated spatial errors do not play any important role  in  the analysis of global errors (space truncation  errors plus time integration errors) as it  can be seen e.g. in \cite[Chapt. IV]{Hun-Ver03}. It should be remarked that the discretization of the source term $c(t,\vec{x})$ is entirely included in  $g_1(t)$ \cite[Sect. 1]{GHH-sinum20}. Besides, the 
terms 
\begin{equation}\label{es-1} \varphi_i(t):=D_i U(t) + g_i(t),\quad t\in[0,t^*],\quad (i=1,\ldots,m),\quad \mbox{\rm satisfy}\quad \sum_{i=1}^m \varphi_i(t)=\dot U(t)\end{equation}
and they are smooth (i.e. they have bounded first and second derivatives independently of the spatial resolution), since  $U(t)=u(t,\vec{x}_G)$ is a smooth function  and we have (below $\delta_{i,j}=0$ for $i\ne j$ and $\delta_{i,i}=1$)
$$  \varphi_i(t):=D_i U(t) + g_i(t)= \beta_i \partial_{x_ix_i} u(t,\vec{x}_G) + \delta_{i1}\cdot c(t,\vec{x}_G) + \mathcal{O}\left((\Delta x_i)^2 \partial_{x_ix_ix_ix_i} u(t,\vec{x}_G)\right).$$  
Additionally, when time independent boundary conditions are assumed in the PDE problem (\ref{lindiff}), we have   $\dot g_i(t)=0,\;i=2,\ldots,m,$ and \cite[Sect. 1 and 4]{GHH-sinum20},
\begin{equation}\label{es-2}
\begin{array}{rl}
D_{i_1} D_{i_2} \ldots D_{i_r} \dot\varphi_j(t)&=  \beta_{i_1} \beta_{i_2} \cdots \beta_{i_r} \dfrac{\partial^{2r+2} \dot u (t,\vec{x}_G)}{\partial {x_{i_1}^2}\partial {x_{i_2}^2} \cdots \partial {x_{i_r}^2}\partial {x_{j}^2}}\\& + \mathcal{O}\left((\Delta x_{i_1})^2 +\cdots +(\Delta x_{i_r})^2+(\Delta x_{j})^2\right), \quad i_1<i_2< \cdots <i_r<j. \\
D_{i_1} D_{i_2} \ldots D_{i_r} \ddot\varphi_j(t)&=  \beta_{i_1} \beta_{i_2} \cdots \beta_{i_r} \dfrac{\partial^{2r+2} \ddot u (t,\vec{x}_G)}{\partial {x_{i_1}^2}\partial {x_{i_2}^2} \cdots \partial {x_{i_r}^2}\partial {x_{j}^2}} \\&+ \mathcal{O}\left((\Delta x_{i_1})^2 +\cdots +(\Delta x_{i_r})^2+(\Delta x_{j})^2\right), \quad i_1<i_2< \cdots <i_r<j.   
\end{array}
\end{equation}

\begin{theorem}\label{5-eq8} Assume that the exact solution of the discretized problem (\ref{ODEm}) satisfies the following uniform bounds
$$ \|U^{(j)}(t)\|_\infty \leq C, \quad \|g^{(i)}(t)\|_\infty \le C, \quad i=0,\ldots,3,\;j=0,\ldots,4,\quad \quad t\in [0,t^*],$$
that $\dot g_j(t)=0, \;j=2,\ldots,m,$ and that  (\ref{es-2}) holds. Then, the global errors, with $D E_0=\mathcal{O}(\tau^2)$, for the AMF-W-method (\ref{AMF-W1-m}) with $\theta=1/2$   fulfill $\| E_n\| \le C' \tau^2,$ $n=1,\ldots,n^*=t^*/\tau,$
 where the constant $C'$ only depends on $m$ and $C$. 
\end{theorem}
{\bf Proof.} 
According to  \cite[formula (2.11)]{GHH-sinum20} the global errors of (\ref{AMF-W1-m}) follow the recursion 
\begin{equation}\label{2.11}
E_n=(R^nD^{-1})D E_0 + \sum_{j=0}^{n-1} R^{n-1-j} S_j, \quad n=1,2,\dots,n^*=t^*/\tau,
\end{equation}
where the matix $R$ is given by (\ref{stab-matrix-m}) and the discretization local errors are given by \cite[formula (2.10)]{GHH-sinum20}
\begin{equation}\label{2.10}
S_n= \Pi(\theta)^{-1} \big(\tau \dot U(t_n) + \theta \tau^2 \dot{\mathcal{G}} (t_n)\big)-\big(U(t_n+\tau)-U(t_n)\big),
\end{equation}
with
%\begin{equation}\label{2.4}
$\mathcal{G}(t)=\sum_{i=1}^m \left(\prod_{j=1}^{i-1} (I-\theta \tau D_j)\right) g_i(t)$, and the convention $\prod_{j=1}^{0} (I-\theta \tau D_j)=I$ (see \cite[formula (2.4)]{GHH-sinum20}).
%\end{equation}
We also make use of other expression for the global errors (see \cite[formula (4.11)]{GHH-sinum20}), obtained by partial summation in (\ref{2.11}),
\begin{equation}\label{4.11}
E_n=(I-R^n)(I-R)^{-1}S_0 + \sum_{j=0}^{n-2} (I-R^{n-1-j})(I-R)^{-1}(S_{j+1}-S_j), \quad n=1,2,\ldots,
\end{equation}
and of a simplified expression for the local  errors given in \cite[formula (4.7)]{GHH-sinum20} 
\begin{equation}\label{4.7}\begin{array}{c}
S_n=S_n^{(1)} + S_n^{(2)},\\[0.5pc]
S_n^{(1)}=\displaystyle{\frac{\tau^2}{2}\Pi^{-1}\sum_{i=1}^m(\Pi_i-\Pi)\dot \varphi_i(t_n)},\quad 
S_n^{(2)}=-\frac{\tau^3}{2}\int_0^1 (1-s)^2\dddot U(t_n+s\tau)ds , \\[0.5pc]
\Pi=(I-\theta \tau D_1)\cdots (I-\theta \tau D_m), \quad \Pi_i:=(I-\theta \tau D_1)\cdots (I-\theta \tau D_{i-1}), \;(i>1), \quad \Pi_1:=I.
\end{array}
\end{equation}
Since, from (\ref{es-1}), $\varphi_1(t)=\dot U(t)-\sum_{i=2}^m \varphi_i(t)$,  
we can split the term $S_n^{(1)}$ of the local error in two parts as
%\begin{equation}\label{5-eq1}
%\varphi_1(t)=\dot U(t)-\sum_{i=2}^m \varphi_i(t),
%\end{equation}
\begin{equation}\label{5-eq2}\begin{array}{c}
S_n^{(1)}=S_n^{(1,a)} + S_n^{(1,b)},\qquad
S_n^{(1,a)}:=\displaystyle\frac{\tau^2}{2}\Pi^{-1}(I-\Pi)\ddot{U}(t_n),\qquad S_n^{(1,b)}:=\frac{\tau^2}{2}\Pi^{-1}\sum_{i=2}^m(\Pi_i-I)\dot \varphi_i(t_n).
\end{array}
\end{equation}
Now, to bound the global errors generated by the contributions of each term of the local error we make use of the following bounds, which are a consequence of Theorem \ref{4-eq3}. In this case, we can apply Theorem \ref{4-eq3} due to the result in \cite{Hun-99} where is guaranted for  $R(z_1,\ldots,z_m)$ that $\alpha=\alpha_m=\frac{\pi}{2(m-1)}, \;m\geq 2$. Consequently it holds that 
\begin{equation}\label{5-eq3}
\| D^{-1}R^n \|_\infty \le C_0, \quad n=0,1,\ldots, n^*=t^*/\tau, \quad \tau>0,\quad N_j \in \mathbb{N}.
\end{equation}
From the assumption $\dot g_j(t)=0, \;j=2,\ldots,m,$ on time-independent boundary conditions (see (\ref{es-2}) and  \cite[Example 4.6]{GHH-sinum20})
\begin{equation}\label{5-eq4}\begin{array}{c}
\|(I-\Pi)\ddot U(t)\|_\infty \le C_1 \tau, \quad  \|(I-\Pi)\dddot U(t)\|_\infty \le C_2 \tau, \\[0.5pc]
\|(I-\Pi_i)\dot \varphi_j(t)\|_\infty \le C_3 \tau, \quad  \|(I-\Pi_i)\ddot \varphi_j(t)\|_\infty \le C_4 \tau, \quad \hbox{\rm when } j\ge i.
\end{array}
\end{equation}
{\sf (1)} 
We start with the global errors $E_n^{(1,a)}$ generated by the local errors  $S_n^{(1,a)}$. To bound them, from (\ref{4.11}) and (\ref{stab-matrix-m})
\begin{equation*}
\begin{split} 
E_n^{(1,a)}&=  \displaystyle (I-R^n)(I-R)^{-1}S_0^{(1,a)} + \sum_{j=0}^{n-2} (I-R^{n-1-j})(I-R)^{-1}(S_{j+1}^{(1,a)}-S_j^{(1,a)})\\&= \displaystyle -\frac{\tau}{2} (I-R^n)D^{-1}(I-\Pi)\ddot U(t_0) -  \frac{\tau}{2} \sum_{j=0}^{n-2} (I-R^{n-1-j})D^{-1}(I-\Pi)\left( \ddot U(t_j+\tau)-\ddot U(t_j)\right)\\&= \mathcal{O}(\tau^2).
\end{split}
\end{equation*}
{\sf (2)} For the global errors  $E_n^{(1,b)}$ generated by the local errors  $S_n^{(1,b)}$, we first take into account that 
\begin{equation}\label{5-eq4a}\begin{array}{c} S_n^{(1,b)}\displaystyle=\frac{\tau^2}{2}\Pi^{-1}\sum_{i=2}^m(\Pi_i-I)\dot \varphi_i(t_n)= \frac{\tau}{2}D^{-1}(R-I)\sum_{i=2}^m(\Pi_i-I)\dot \varphi_i(t_n).\end{array}
\end{equation}
Then, it follows from (\ref{4.11}) and (\ref{stab-matrix-m}) that 
\begin{equation*}
\begin{split}
E_n^{(1,b)}&=  \displaystyle (I-R^n)(I-R)^{-1}S_0^{(1,b)} + \sum_{j=0}^{n-2} (I-R^{n-1-j})(I-R)^{-1}(S_{j+1}^{(1,b)}-S_j^{(1,b)})\\
&=\displaystyle -\frac{\tau}{2} (I-R^n)D^{-1}\sum_{i=2}^m (\Pi_i-I)\dot \varphi_i(t_0)   -  \frac{\tau}{2} \sum_{i=2}^m \left(\sum_{j=0}^{n-2} (I-R^{n-1-j})D^{-1}(\Pi_i-I)\left( \dot \varphi_i(t_{j}+\tau)-\dot \varphi_i(t_j)\right) \right)\\&= \mathcal{O}(\tau^2).
\end{split}
\end{equation*}
{\sf (3)} For the global errors  $E_n^{(2)}$ generated by the local errors  $S_n^{(2)}$, we use the formula (\ref{2.11}). In this case we define 
\begin{equation}\label{5-eq5}\begin{array}{c} V(t):= D\left(-\frac{1}{2} \int_0^1 (1-s)^2\dddot {U}(t +s \tau) ds\right).\end{array}
\end{equation}
so that  
%\begin{equation}\label{5-eq5a}\begin{array}{c} 
$S_n^{(2)}=\tau^3 D^{-1} V(t_n).$
%\end{array}\end{equation}
Besides, 
\begin{equation}\label{5-eq6}
\begin{array}{c} V(t)= -\frac{1}{2} \int_0^1 (1-s)^2 D\dddot{U}(t +s \tau) ds=  -\frac{1}{2} \int_0^1 (1-s)^2\left(U^{(4)}(t +s \tau)-\dddot g(t+s\tau)\right) ds. \end{array}
\end{equation}
From here we deduce (under the regularity assumption in the exact solution) that 
% \begin{equation}\label{5-eq7}\begin{array}{c}  
$\|V(t)\|_\infty \le C_5.$
%\end{array} \end{equation}
Now, the bound for the global errors follows from
\begin{equation} E_n^{(2)}= \sum_{j=0}^{n-1} R^{n-1-j} S_j^{(2)} =\tau^3\sum_{j=0}^{n-1} R^{n-1-j} D^{-1} V(t_j)= \mathcal{O}(\tau^2). \end{equation} \hfill $\Box$

\begin{remark}{\rm
Second order of convergence in the maximum norm
for the Douglas method with $\theta=\frac{1}{2}$  is proved in \cite[Theorem 3.1]{arraras17mds}  under the assumption of power-boundedness for the stability matrix $R$ in (\ref{stab-matrix-m}) and assuming that  \cite[(3.16b), p.\,271]{arraras17mds}
\begin{equation}\label{assump-316b}
\tau^{k-1}D^{-1}  D_{l_1}D_{l_2}\dotsb D_{l_k} v(t_n)=\mathcal{O}(1) ,
\quad 1\le l_1<\ldots < l_k<i \le m ,\quad (v=\dot \varphi_i, \ddot\varphi_i).
\end{equation}
 This assumption was also useful in \cite[Theorem 3.2]{hundsdorfer18oms} in order to prove convergence for linear multistep methods with stabilizing corrections applied to split ODEs. Although (\ref{assump-316b}) is closely related to (\ref{es-2}), the  proof of convergence presented in \cite{arraras17mds}  does  require the assumption on power-boundedness for the stability matrix $R$, which, as far as we are aware,  has not been shown for $m\ge 3$  so far, whereas our proof does not require such an assumption. \hfill $\Box$
}
\end{remark}

\noindent Convergence of order two in the maximum norm  for the Douglas method \cite[p. 373] {Hun-Ver03} and time independent boundary conditions can also be shown following similar steps as in the proof of Theorem \ref{5-eq8}. To this aim, let us consider the Douglas method applied to (\ref{ODEm}):

\begin{equation}\label{Douglas}
\begin{split}
v_0&= U_n+\tau (D U_n+g(t_n))\\
v_i&= v_{i-1}+\theta\tau ((D_i v_i+g_i(t_{n+1}))-(D_i U_n+g_i(t_{n}))),\quad i=1,\ldots,m,\\
U_{n+1}&=v_m.
\end{split}
\end{equation}

\begin{theorem}\label{thm-dou} Under the same assumptions of Theorem \ref{5-eq8}, the global errors, with $D E_0=\mathcal{O}(\tau^2)$, for the Douglas method (\ref{Douglas}) with $\theta=1/2$   fulfill $\| E_n\| \le C' \tau^2,$ $n=1,\ldots,n^*=t^*/\tau,$
 where the constant $C'$ only depends on $m$ and $C$. 
\end{theorem}
{\bf Proof.} The global errors $E_n=U_n-U(t_n)$ for the method (\ref{Douglas}) fulfill the recursion $E_{n+1}=RE_n+S_n$, $n\geq 0$, where the stability matrix $R$ is given by (\ref{stab-matrix-m}) and the local errors $S_n$ are obtained as given in \cite[formula (3.15)]{Hun-Ver03} by 

\begin{equation}\label{locerr-dou}
S_n=-Q_m^{-1}\cdots Q_1^{-1}(r_0+r_1)-Q_m^{-1}\cdots Q_2^{-1}r_2-\ldots-Q_m^{-1}r_m,
\end{equation}
with $Q_i=I-\theta\tau D_i$, $1\leq i\leq m$, whereas, taking into account that in (\ref{Douglas}) $v_i\approx U(t_n+1)$, $0\leq i\leq m$, for $r_i=r_i(t_n)$ it holds that
\begin{equation}\label{r0}
\begin{split}
r_0(t_n)&=U(t_{n+1})-U(t_n)-\tau \dot U(t_n)=\frac{\tau^2}{2}\ddot U(t_n)+\frac{\tau^3}{2}\int_0^1 (1-s)^2\dddot U(t_n+s\tau)ds,
\end{split}
\end{equation}
and, for $1\leq i\leq m$, 
\begin{equation}\label{ri}
\begin{split}
r_i(t_n)&=-\theta\tau\left( \varphi_i(t_{n+1})- \varphi_i(t_{n})\right)= -\theta\tau^2 \int_0^1 \dot \varphi_i(t_n+s\tau)ds=-\theta\tau^2 \dot\varphi_i(t_n)-\theta\tau^3\int_0^1 (1-s) \ddot \varphi_i(t_n+s\tau)ds.
\end{split}
\end{equation}
With $\Pi$ and $\Pi_i$ defined in (\ref{4.7}), we can rewrite (\ref{locerr-dou}) as
\begin{equation}\label{locerr-dou-1}
\begin{split}
S_n&=-\Pi^{-1}\left((r_0+r_1)+Q_1r_2+Q_1Q_2 r_3+\ldots+Q_1\cdots Q_{m-1} r_m\right)\\
&=-\Pi^{-1}\left(\big(\sum_{i=0}^m r_i\big)+(\Pi_2-I)r_2+(\Pi_3-I) r_3+\ldots+(\Pi_m-I) r_m\right),\quad r_i=r_i(t_n),\;0\leq i\leq m.
\end{split}
\end{equation}
With $\theta=\frac{1}{2}$, using (\ref{r0}), (\ref{ri}) and (\ref{es-1}), it holds that 
\begin{equation}\label{r0-1}
\begin{split}
\sum_{i=0}^m r_i(t_n)&=\frac{\tau^3}{2}\int_0^1 (s^2-s)\dddot U(t_n+s\tau)ds.
\end{split}
\end{equation}
Now, we split the local error in two parts as 
\begin{equation}\label{locerr-dou-split}\begin{array}{c}
S_n=S_n^{(a)} + S_n^{(b)},\\
S_n^{(a)}:=-\Pi^{-1}\big(\sum_{i=0}^m r_i\big),\quad S_n^{(b)}:=-\Pi^{-1}\left((\Pi_2-I)r_2+(\Pi_3-I) r_3+\ldots+(\Pi_m-I) r_m\right).
\end{array}
\end{equation}
Partial summation in the global error recursion leads us to the relation (\ref{4.11}), with $(I-R)^{-1}\Pi^{-1}=-\tau^{-1}D^{-1}.$ Now, we bound the global errors generated by the contributions of each term of the local error.

\medskip

\noindent {\sf (a)} 
For the global errors $E_n^{(a)}$ generated by the local errors  $S_n^{(a)}$, using (\ref{4.11}), (\ref{r0-1}) and (\ref{locerr-dou-split}) we have that

\begin{equation*}
\begin{split}
E_n^{(a)}&= \displaystyle (I-R^n)D^{-1}\tau^{-1} \big(\sum_{i=0}^m r_i(t_0)\big)+ \sum_{j=0}^{n-2} (I-R^{n-1-j})D^{-1}\tau^{-1} \big(\sum_{i=0}^m r_i(t_{j+1})-r_i(t_{j})\big)\\&= \displaystyle \frac{\tau^2}{2} (I-R^n)D^{-1}\int_0^1 (s^2-s)\dddot U(t_0+s\tau)ds\\&+\frac{\tau^2}{2} \sum_{j=0}^{n-2} (I-R^{n-1-j})D^{-1}\int_0^1 (s^2-s)\big(\dddot U(t_{j+1}+s\tau)-\dddot U(t_j+s\tau)\big)ds\\&= \mathcal{O}(\tau^2).
\end{split}
\end{equation*}

\noindent {\sf (b)} 
For the global errors $E_n^{(b)}$ generated by the local errors  $S_n^{(b)}$, using (\ref{4.11}), (\ref{ri}), (\ref{locerr-dou-split}) and (\ref{5-eq4}) it holds that

\begin{equation*}
\begin{split}
E_n^{(b)}&= \displaystyle (I-R^n)D^{-1}\tau^{-1} \sum_{i=2}^m (\Pi_i-I) r_i(t_0)+ \sum_{j=0}^{n-2} (I-R^{n-1-j})D^{-1}\tau^{-1} \sum_{i=2}^m (\Pi_i-I) (r_i(t_{j+1})-r_i(t_{j}))\\&= \displaystyle \frac{\tau}{2} (I-R^n)D^{-1}\sum_{i=2}^m \int_0^1 (I-\Pi_i)\dot \varphi_i(t_0+s\tau)ds\\&+\frac{\tau}{2} \sum_{i=2}^m \sum_{j=0}^{n-2} (I-R^{n-1-j})D^{-1}\int_0^1 (I-\Pi_i)\big(\dot \varphi_i(t_{j+1}+s\tau)-\dot \varphi_i(t_{j}+s\tau)\big)ds\\&= \mathcal{O}(\tau^2).
\end{split}
\end{equation*}
\hfill $\Box$

\begin{remark}{\rm
A modified one-stage AMF-W method (henceforth denoted as {\bf modified AMF-W1})
\begin{equation}\label{AMF-W1-mod-m}
\begin{aligned}
K_1^{(0)} &= \displaystyle \tau D\,  U_n 
+ \tau\, g(t_n)   \\[1.5mm]
(I-\theta \tau D_j ) K_1^{(j)} &= K_1^{(j-1)} + \theta \tau^2 \dot g_j (t_n+\tau/2),\quad j=1,\ldots,m,\\[1.5mm]
\displaystyle U_{n+1} &= U_n + K_1^{(m)} ,
\end{aligned}
\end{equation}
was introduced in \cite{GHH-sinum20}. For this method, second order convergence in the $\ell_\infty-$norm for time independent boundary conditions can be shown as in Theorem \ref{5-eq8} above considering that its local error can be expressed as (see \cite[formula (5.4)]{GHH-sinum20} )

\begin{equation*}
\begin{split}
S_n  &= \tau^3 \Pi^{-1}\int_0^{1/2} \!\! s\, \dddot U(t_n + s\tau ) \, \d s +  \displaystyle \frac{\tau^2}2 
\Pi^{-1}\Bigl( \,\sum_{i=2}^m \Bigl( \Pi_i  -I\Bigr) \,
 \dot \varphi_i (t_n + \tfrac \tau 2 ) 
\Bigr)   \\  &-  \displaystyle  \frac{\tau^3}{2}\int_0^1 k(s)\, \dddot U(t_n + s\tau ) \, \d s   , \quad k(s) := \min \bigl( s^2, (1-s)^2 \bigr).
\end{split}
\end{equation*} \hfill $\Box$

}
\end{remark}

\begin{theorem}\label{thm-modAMF-W1} Under the same assumptions of Theorem \ref{5-eq8}, the global errors, with $D E_0=\mathcal{O}(\tau^2)$, for the modified AMF-W method (\ref{AMF-W1-mod-m}) with $\theta=1/2$   fulfill $\| E_n\| \le C' \tau^2,$ $n=1,\ldots,n^*=t^*/\tau,$
 where the constant $C'$ only depends on $m$ and $C$. 
\end{theorem}

\noindent {\bf Proof.} The proof follows along the lines of the proofs of Theorems \ref{5-eq8} and \ref{thm-dou}  \hfill $\Box$

\begin{remark}{\rm
For $m\geq 3$ spatial dimensions, the three methods (\ref{AMF-W1-m}), (\ref{AMF-W1-mod-m}) and (\ref{Douglas}) display first order of convergence in the $\ell_\infty$-norm (up to a logarithmic factor) when applied to (\ref{ODEm}) with time dependent Dirichlet boundary conditions. When $m=2$, the methods (\ref{AMF-W1-mod-m}) and (\ref{Douglas}) have the advantage that their global error in the $\ell_\infty$-norm is $\min\{\mathcal{O}(\tau^2|\log\;h|^2),\: \mathcal{O}(\tau)\}$ (see \cite[Section 5]{GHH-sinum20}), whereas the method in (\ref{AMF-W1-m})  only has convergence of size $\mathcal{O}(\tau)$ in the maximum norm . \hfill $\Box$
}
\end{remark}

\section{Numerical Illustration}

\noindent We first consider the linear diffusion partial differential equation (\ref{lindiff}) in three and four spatial dimensions, with diffusion coefficients $\beta_j=1$, $1\leq j\leq m$. Our aim is to  illustrate numerically the second order convergence in the maximum norm for the one-stage AMF-W method (\ref{AMF-W1-m}) (and its modified version (\ref{AMF-W1-mod-m})) and the Douglas method (\ref{Douglas}), both with parameter $\theta=\frac{1}{2}$, when time independent boundary conditions are imposed on the PDE. For time dependent boundary conditions, order one (up to a logarithmic factor) is attained by both methods when the spatial dimension is $m\geq 3$. For our numerical experiments we consider that $c(t,\vec{x})$ is selected in such way that
\begin{equation}\label{exactsol}
%\textstyle
u(t,\vec{x})=u_e(t,\vec{x}):= \e^t \bigg( 4^m \prod_{j=1}^m x_j(1-x_j) + \kappa \sum_{j=1}^m \Bigl(x_j+\frac{1}{j+2}\Bigr)^2\bigg)
\end{equation}
is the exact solution of (\ref{lindiff}). We impose the initial condition $u(0,\vec{x})=u_e(0,\vec{x})$
and Dirichlet boundary conditions. Here, we consider the cases $m=3,4$. If $\kappa=0$ we have homogeneous boundary conditions,
but when $\kappa=1$ we get non-homogeneous time-dependent Dirichlet conditions.

\medskip

\begin{figure}[h!]
\epsfig{figure=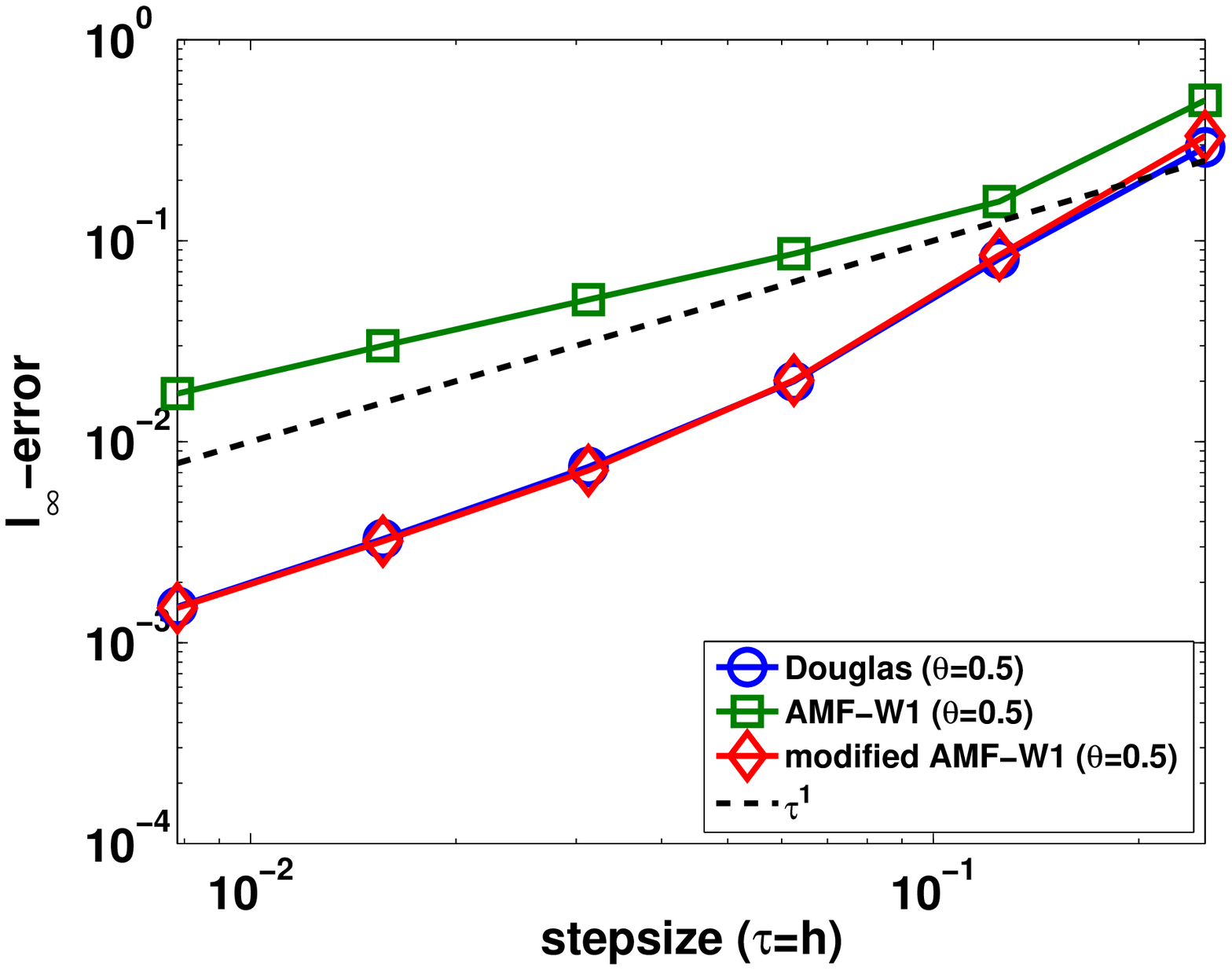,scale=0.43}\hspace{0.5cm}
\epsfig{figure=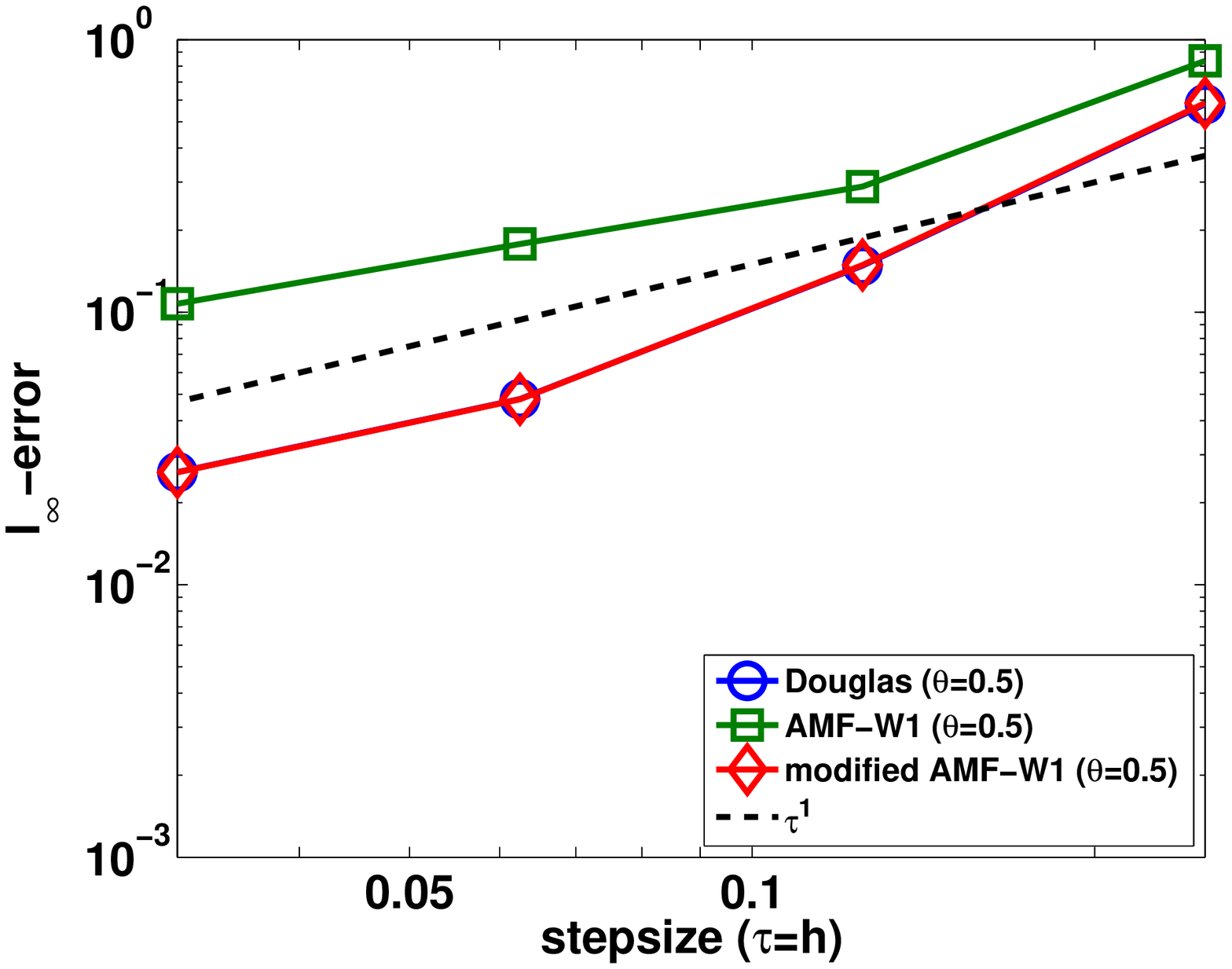,scale=0.43}
\caption{Error in the $\ell_\infty-$norm vs stepsize on the linear model (\ref{lindiff})-(\ref{exactsol}) with time dependent boundary conditions ($\kappa=1$) and $\tau=h=\Delta x_i$, $1\leq i\leq m$.
Spatial dimension $m=3$ (left) and $m=4$ (right). A dashed straight line with slope one is included to compare the PDE order of convergence.}
\label{fig-3d4dk1}
\end{figure}

\noindent  We apply the MOL approach on a uniform grid with meshwidth $h=\Delta x_i=1/(N+1)$, $1\leq i\leq m$, where $N=2^j-1$, $j=2,\ldots,j_{{\rm max}}$, with $j_{{\rm max}}=7$ for $m=3$ and $j_{{\rm max}}=5$ if $m=4$. Hence, a semi-discretized system with corresponding dimension $N^m$ of the form (\ref{ODEm})
is obtained, where $D$ is given in (\ref{4-eq1}) and $g(t)$ includes the discretization of the term $c(t,\vec{x})$ and the terms due to non-homogeneous boundary conditions. Observe that the exact solution (\ref{exactsol}) is a polynomial
of degree $2$ in each spatial variable so that the global errors  come only  from the time discretization. The methods (\ref{AMF-W1-m}), (\ref{AMF-W1-mod-m}) and (\ref{Douglas}) are then applied to (\ref{ODEm}) with fixed step size $\tau=h=2^{-j}$, $2\leq j\leq j_{{\rm max}}$, and the corresponding global errors regarding the PDE solution versus the stepsize are displayed below in Figure \ref{fig-3d4dk0} in case of time independent boundary conditions ($\kappa=0$) and in Figure \ref{fig-3d4dk1} in case of time dependent boundary conditions ($\kappa=1$). In the first case, all methods display second order convergence in the $\ell_\infty$-norm in both dimensions $m=3$ and $m=4$. In the second situation with time dependent boundary conditions, all methods suffer an order reduction and the corresponding  orders of convergence are at most one.

\begin{figure}[h!]
\epsfig{figure=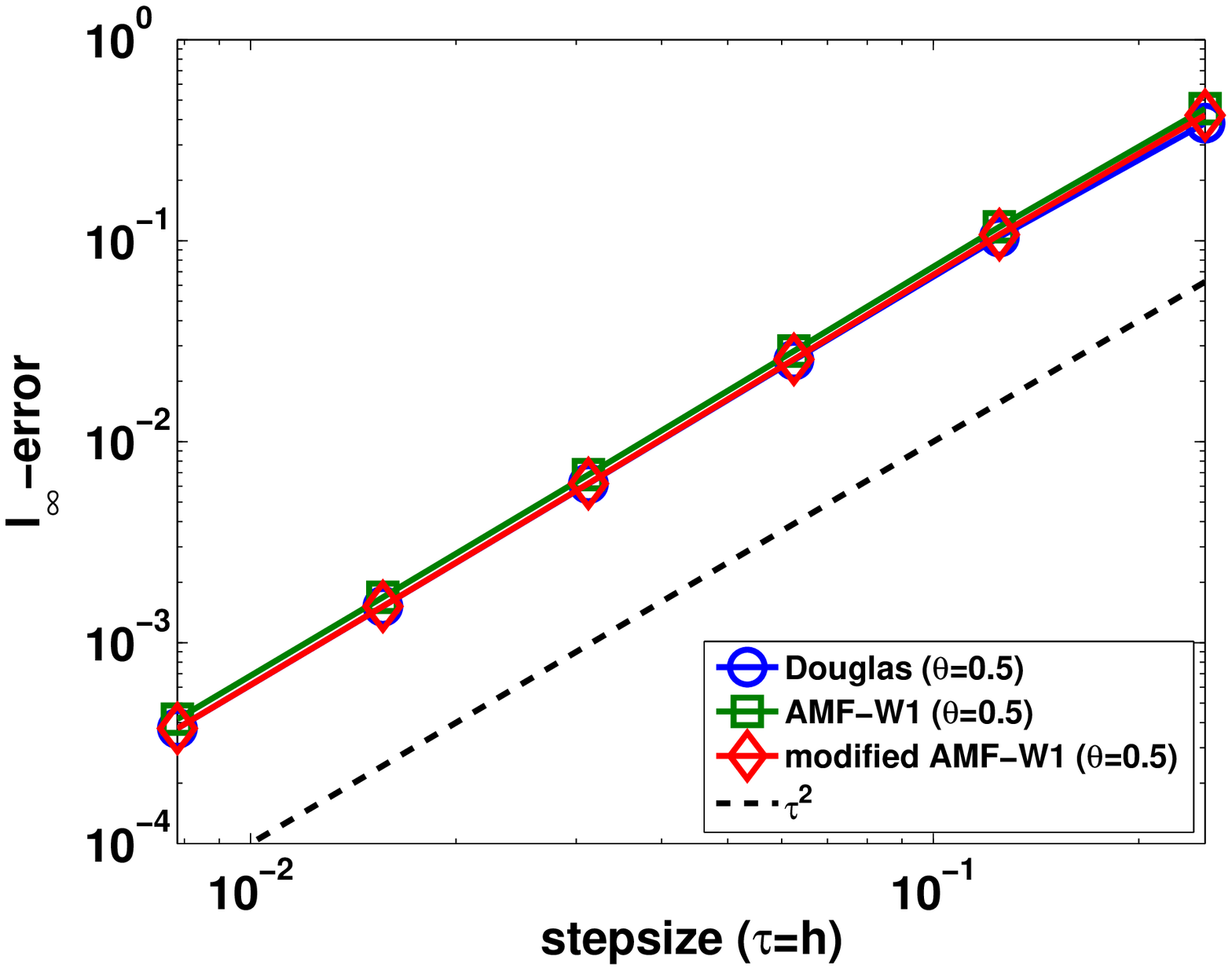,scale=0.43}\hspace{0.5cm}
\epsfig{figure=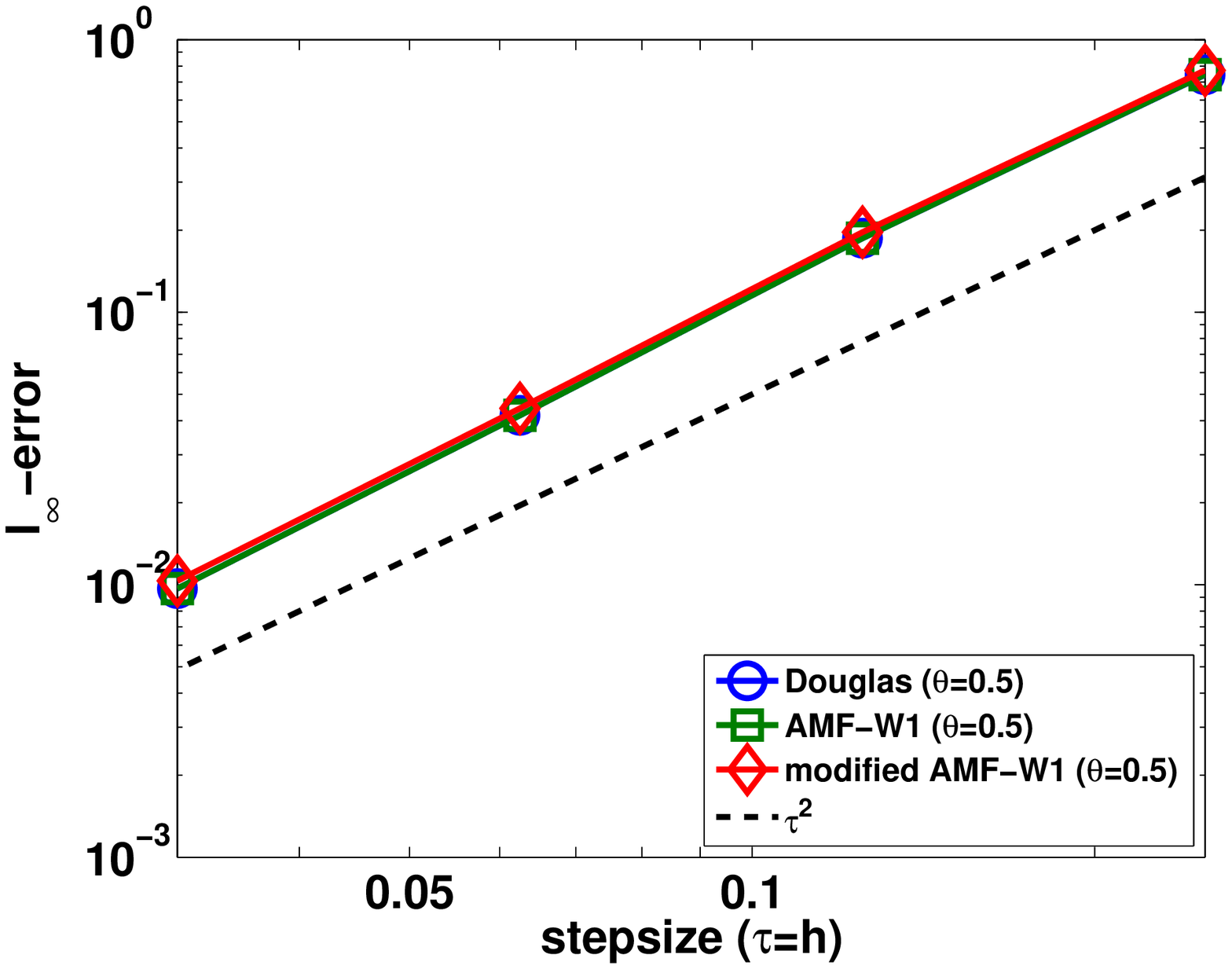,scale=0.43}
\caption{Error in the $\ell_\infty-$norm vs stepsize on the linear model (\ref{lindiff})-(\ref{exactsol}) with time independent boundary conditions ($\kappa=0$) and $\tau=h=\Delta x_i$, $1\leq i\leq m$.
Spatial dimension $m=3$ (left) and $m=4$ (right). A dashed straight line with slope two is included to compare the PDE order of convergence.}
\label{fig-3d4dk0}
\end{figure}

\medskip

A second numerical experiment is included below in Figure \ref{fig-3dvark0k1} for the case of variable diffusion coefficients $\beta_j=\beta_j(\vec{x})$, $1\leq j\leq m$. Although a theoretical analysis for such a case lies beyond the scope of this paper, similar orders of convergence are observed. To illustrate this assertion, we consider $m=3$ spatial dimensions and diffusion coefficients
\begin{equation}\label{vardifcoef}
\beta_1=\beta_1(x,y,z)=(1+xyz)^2,\quad \beta_2=\beta_2(x,y,z)=e^{x-2y+3z} \quad {\rm and} \quad \beta_3=\beta_3(x,y,z)=(1+x^2)e^{-y^2 z}.
\end{equation}
Again $c(t,\vec{x})$ is selected in such way that (\ref{exactsol}) is the exact solution of (\ref{lindiff}), with  homogeneous boundary conditions if $\kappa=0$ 
and time-dependent Dirichlet conditions when $\kappa=1$. The MOL approach is applied on a uniform grid with meshwidth $h=\Delta x_i=1/(N+1)$, $1\leq i\leq 3$, where $N=2^j-1$, $j=2,\ldots,7$. The results displayed in Figure \ref{fig-3dvark0k1} (left) show that all methods provide second order of convergence in the $\ell_\infty$ norm when time-independent boundary conditions are considered ($\kappa=0$). For the case of time-dependent boundary conditions ($\kappa=1$), Figure \ref{fig-3dvark0k1} (right) show an order reduction to at most order one for the three methods considered.

\begin{figure}[h!]
\epsfig{figure=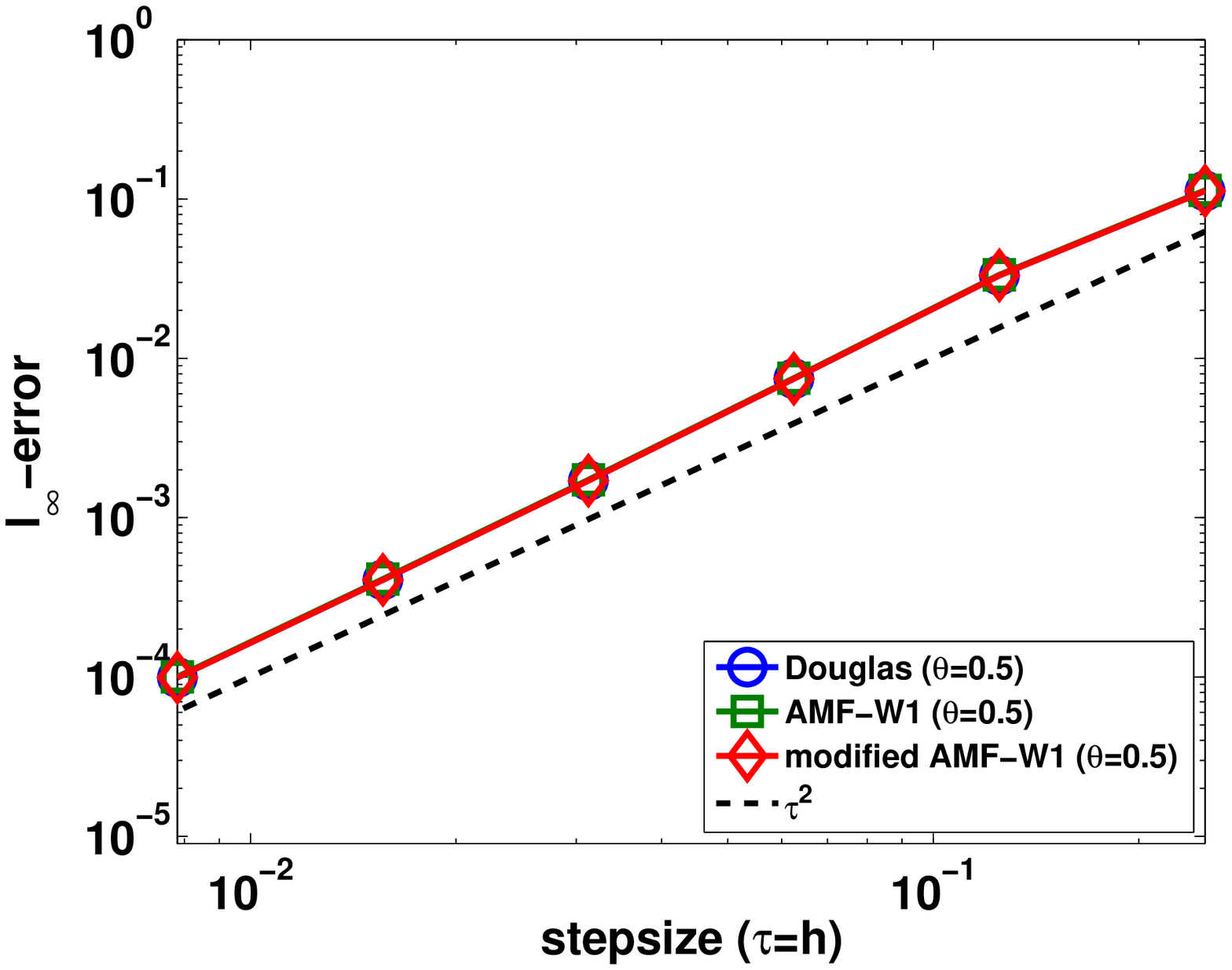,scale=0.43}\hspace{0.5cm}
\epsfig{figure=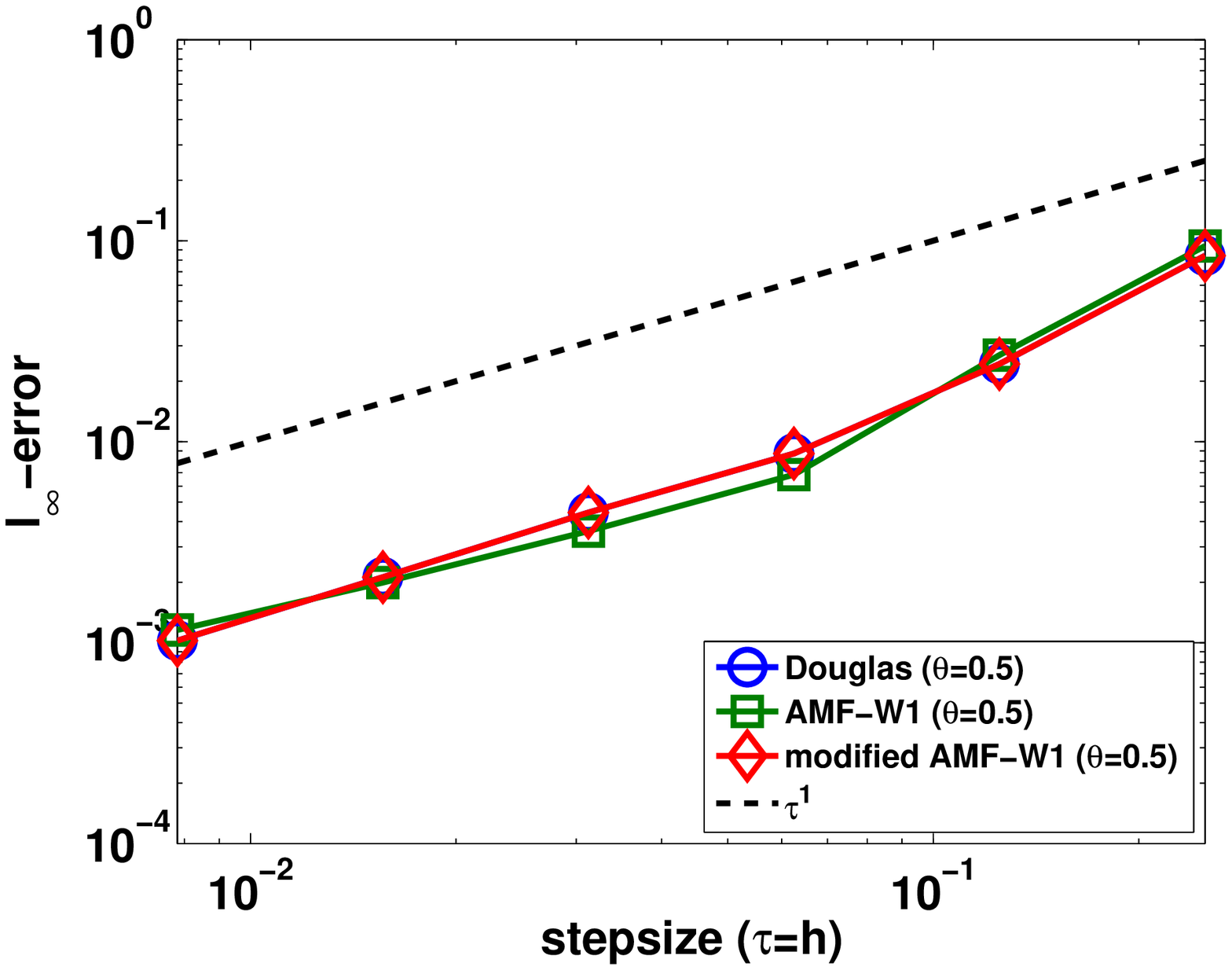,scale=0.43}
\caption{Error in the $\ell_\infty-$norm vs stepsize on the linear model (\ref{lindiff})-(\ref{exactsol}) with variable diffusion coefficients (\ref{vardifcoef}), $\tau=h=\Delta x_i$, $1\leq i\leq m$, and spatial dimension $m=3$. Time independent boundary conditions $\kappa=0$ (left) and time dependent boundary conditions $\kappa=1$ (right). Dashed straight lines with slopes two and one, respectively, are included to compare the PDE order of convergence.}
\label{fig-3dvark0k1}
\end{figure}

\end{document}